\documentclass{amsart}

\usepackage{amssymb,mathscinet}
\usepackage{amsthm}
\usepackage{amsmath}
\usepackage[all]{xy} \SelectTips{eu}{}
\usepackage[hidelinks]{hyperref}
\newcommand{\numberseries}{\bfseries}   %Fontseries used for numbering

\newlength{\thmtopspace}                %Space above theorem
\newlength{\thmbotspace}                %Space below theorem
\newlength{\thmheadspace}               %Space after theorem label
\newlength{\thmindent}                  %For indenting

\newtheoremstyle{bfupright head,upright body}
                {\thmtopspace}{\thmbotspace}
                {\upshape}{\thmindent}{\bfseries}{.}{\thmheadspace}
                {{\numberseries \thmnumber{#2\;}}\thmnote{#3}}

\newtheoremstyle{fixed bf head,slanted body}
                {\thmtopspace}{\thmbotspace}{\slshape}
                {\thmindent}{\bfseries}{.}{\thmheadspace}
                {{\numberseries \thmnumber{#2\;}}\thmname{#1}\thmnote{ (#3)}}

\newtheoremstyle{fixed bf head,upright body}
                {\thmtopspace}{\thmbotspace}{\upshape}
                {\thmindent}{\bfseries}{.}{\thmheadspace}
                {{\numberseries \thmnumber{#2\;}}\thmname{#1}\thmnote{ (#3)}}

\theoremstyle{bfupright head,upright body}
\newtheorem{res}{}[section]
\newtheorem{bfhpg}[res]{}

\theoremstyle{fixed bf head,slanted body}
\newtheorem{thm}[res]{Theorem}          \newtheorem*{thm*}{Theorem}
\newtheorem{prp}[res]{Proposition}      \newtheorem*{prp*}{Proposition}
\newtheorem{cor}[res]{Corollary}        \newtheorem*{cor*}{Corollary}
\newtheorem{lem}[res]{Lemma}            \newtheorem*{lem*}{Lemma}

\theoremstyle{fixed bf head,upright body}
\newtheorem{dfn}[res]{Definition}       \newtheorem*{dfn*}{Definition}
\newtheorem{rmk}[res]{Remark}           \newtheorem*{rmk*}{Remark}
\newtheorem{exa}[res]{Example}          \newtheorem*{exa*}{Example}
\newtheorem{fct}[res]{Fact}             \newtheorem*{fct*}{Fact}

\newlength{\thmlistleft}        %leftmargin
\newlength{\thmlistright}       %rightmargin
\newlength{\thmlistpartopsep}   %partopsep
\newlength{\thmlisttopsep}      %topsep
\newlength{\thmlistparsep}      %parsep
\newlength{\thmlistitemsep}     %itemsep

\newcounter{eqc}
\newenvironment{eqc}{\begin{list}{\upshape (\textit{\roman{eqc}})}%
    {\usecounter{eqc}%
      \setlength{\leftmargin}{\thmlistleft}%
      \setlength{\labelwidth}{\thmlistleft}%
      \setlength{\rightmargin}{\thmlistright}%
      \setlength{\partopsep}{\thmlistpartopsep}%
      \setlength{\topsep}{\thmlisttopsep}%
      \setlength{\parsep}{\thmlistparsep}%
      \setlength{\itemsep}{\thmlistitemsep}}}%
  {\end{list}}%

\newcommand{\eqclbl}[1]{{\upshape(\textit{#1})}}

\newcounter{prt}
\newenvironment{prt}{\begin{list}{\upshape (\alph{prt})}%
    {\usecounter{prt}%
      \setlength{\leftmargin}{\thmlistleft}%
      \setlength{\labelwidth}{\thmlistleft}%
      \setlength{\rightmargin}{\thmlistright}%
      \setlength{\partopsep}{\thmlistpartopsep}%
      \setlength{\topsep}{\thmlisttopsep}%
      \setlength{\parsep}{\thmlistparsep}%
      \setlength{\itemsep}{\thmlistitemsep}}}%
  {\end{list}}%

\newcounter{rqm}
\newenvironment{rqm}{\begin{list}{\upshape (\arabic{rqm})}%
    {\usecounter{rqm}%
      \setlength{\leftmargin}{\thmlistleft}%
      \setlength{\labelwidth}{\thmlistleft}%
      \setlength{\rightmargin}{\thmlistright}%
      \setlength{\partopsep}{\thmlistpartopsep}%
      \setlength{\topsep}{\thmlisttopsep}%
      \setlength{\parsep}{\thmlistparsep}%
      \setlength{\itemsep}{\thmlistitemsep}}}%
  {\end{list}}%

\newenvironment{itemlist}{\nopagebreak \begin{list}{$\bullet$}%
    {\setlength{\leftmargin}{\thmlistleft}%
      \setlength{\labelwidth}{\thmlistleft}%
      \setlength{\rightmargin}{\thmlistright}%
      \setlength{\partopsep}{\thmlistpartopsep}%
      \setlength{\topsep}{\thmlisttopsep}%
      \setlength{\parsep}{\thmlistparsep}%
      \setlength{\itemsep}{\thmlistitemsep}}}%
  {\end{list}}%

\newenvironment{prf*}[1][Proof]{%
  \begin{proof}[\bf #1]
    \setcounter{equation}{0}
    }
  {\end{proof}
}

\newcommand{\proofofimp}[3][:]{\mbox{\eqclbl{#2}$\!\implies\!$\eqclbl{#3}#1}}

\newcommand{\pgref}[1]{\ref{#1}}
\newcommand{\thmref}[2][Theorem~]{#1\pgref{thm:#2}}
\newcommand{\corref}[2][Corollary~]{#1\pgref{cor:#2}}
\newcommand{\prpref}[2][Proposition~]{#1\pgref{prp:#2}}
\newcommand{\lemref}[2][Lemma~]{#1\pgref{lem:#2}}
\newcommand{\dfnref}[2][Definition~]{#1\pgref{dfn:#2}}
\newcommand{\exaref}[2][Example~]{#1\pgref{exa:#2}}
\newcommand{\rmkref}[2][Remark~]{#1\pgref{rmk:#2}}
\newcommand{\secref}[2][Section~]{#1\ref{sec:#2}}

\newcommand{\thmcite}[2][?]{\cite[thm.~#1]{#2}}
\newcommand{\corcite}[2][?]{\cite[cor.~#1]{#2}}
\newcommand{\prpcite}[2][?]{\cite[prop.~#1]{#2}}
\newcommand{\lemcite}[2][?]{\cite[lem.~#1]{#2}}

\newcommand{\dfncite}[2][?]{\cite[def.~#1]{#2}}
\renewcommand{\eqref}[1]{(\pgref{eq:#1})}
\numberwithin{equation}{res}

\def\urltilda{\kern -.15em\lower .7ex\hbox{\~{}}\kern .04em}
\newcommand{\Cone}[1]{\nobreak{\operatorname{Cone}#1}}
\newcommand{\cathom}[3]{\operatorname{hom}_{#1}(#2,#3)}
\newcommand{\setof}[3][\mspace{1mu}]{\{#1#2 \mid #3#1\}}
\newcommand{\Ext}[4][A]{\operatorname{Ext}_{#1}^{#2}(#3,#4)}

\newcommand{\NN}{\mathbb{N}}
\newcommand{\ZZ}{\mathbb{Z}}
\newcommand{\QQ}{\mathbb{Q}}
\newcommand{\Cy}[2][]{\operatorname{Z}_{#1}(#2)}

\newcommand{\lra}{\longrightarrow}
\newcommand{\xla}[2][]{\xleftarrow[#1]{\;#2\;}}
\newcommand{\xra}[2][]{\xrightarrow[#1]{\;#2\;}}
\newcommand{\deq}{\:=\:}

\newcommand{\qqtext}[1]{\qquad\text{#1}\qquad}

\newcommand{\qqand}{\qqtext{and}}
\newcommand{\gra}{\alpha}
\newcommand{\gre}{\varepsilon}
\newcommand{\gri}{\iota}
\newcommand{\grs}{\sigma}
\newcommand{\grz}{\zeta}

\newcommand{\grf}{\varphi}
\newcommand{\mapdef}[4][\rightarrow]{\nobreak{#2\colon #3 #1 #4}}
\newcommand{\dmapdef}[4][\lra]{\nobreak{#2\colon #3\:#1\:#4}}
\newcommand{\is}{\cong}
\newcommand{\qis}{\simeq}

\newcommand{\Aop}{A^\circ}
\newcommand{\sfA}{\mathsf{A}}

\newcommand{\sfU}{\mathsf{U}}
\newcommand{\sfV}{\mathsf{V}}
\newcommand{\sfW}{\mathsf{W}}

\newcommand{\sfUperp}{\mathsf{U}^\perp}
\newcommand{\sfperpV}{{}^\perp\mathsf{V}}
\newcommand{\sfUV}{(\sfU,\sfV)}
\newcommand{\capUV}{\sfU\cap\sfV}

\newcommand{\capUU}{\sfU\cap\sfU^\perp}
\newcommand{\capVV}{{}^\perp\sfV\cap\sfV}

\newcommand{\tp}[3][A]{\nobreak{#2\otimes_{#1}#3}}

\newcommand{\Hom}[3][A]{\operatorname{Hom}_{#1}(#2,#3)}

\newcommand{\fuTdot}{\operatorname{\dot{T}_R}}
\newcommand{\fuT}{\operatorname{T_R}}

\newcommand{\fuI}{\operatorname{I}}

\newcommand{\catGL}[2][\sfV]{\mathsf{LGor}_{#1}(#2)}
\newcommand{\catGR}[2][\sfU]{\mathsf{RGor}_{#1}(#2)}
\newcommand{\catG}[2][\sfW]{\mathsf{Gor}_{#1}(#2)}
\newcommand{\catGFlat}[1]{\mathsf{GFlat}(#1)}
\newcommand{\catCot}[1]{\mathsf{Cot}(#1)}
\newcommand{\catFlat}[1]{\mathsf{Flat}(#1)}
\newcommand{\catFlatCot}[1]{\mathsf{FlatCot}(#1)}
\newcommand{\catPrj}[1]{\mathsf{Prj}(#1)}
\newcommand{\catInj}[1]{\mathsf{Inj}(#1)}
\newcommand{\catMod}[1]{\mathsf{Mod}(#1)}
\newcommand{\stabcatGR}[2][\sfU]{\mathsf{StRGor}_{#1}(#2)}
\newcommand{\stabcatGL}[2][\sfV]{\mathsf{StLGor}_{#1}(#2)}
\newcommand{\stabcatG}[2][\sfW]{\mathsf{StGor}_{#1}(#2)}
\newcommand{\sfCot}{\mathsf{Cot}}
\newcommand{\sfFlat}{\mathsf{Flat}}
\newcommand{\sfFlatCot}{\mathsf{FlatCot}}

\newcommand{\catK}[1]{\mathsf{K}(#1)}
\newcommand{\catKFtac}[1]{\mathsf{K}_{\textnormal{F-tac}}(#1)}
\newcommand{\catKRtac}[2][\sfU]{\mathsf{K}^\textnormal{R}_{#1\textnormal{-tac}}(#2)}
\newcommand{\catKLtac}[2][\sfV]{\mathsf{K}^\textnormal{L}_{#1\textnormal{-tac}}(#2)}
\newcommand{\catKpac}[1]{\mathsf{K}_{\textnormal{pac}}(#1)}
\newcommand{\catKtac}[1]{\mathsf{K}_{\textnormal{tac}}(#1)}
\newcommand{\Ftac}{\textnormal{\bf F}-totally acyclic}
\newcommand{\scR}{\textsc{r}}
\newcommand{\scL}{\textsc{l}}

\title[Homotopy categories of totally acyclic complexes]{Homotopy
  categories of totally acyclic complexes with applications to
  the flat--cotorsion theory}

\author[L.W.\ Christensen]{Lars Winther Christensen}

\address{L.W.C. \ Texas Tech University, Lubbock, TX 79409, U.S.A.}

\email{lars.w.christensen@ttu.edu}

\urladdr{http://www.math.ttu.edu/\urltilda lchriste}

\author[S.\ Estrada]{Sergio Estrada}

\address{S.E. \ Universidad de Murcia, Murcia 30100, Spain}

\email{sestrada@um.es}

\urladdr{http://webs.um.es/sestrada}

\author[P.\ Thompson]{Peder Thompson}

\address{P.T. \ Norwegian University of Science and Technology, 7491
  Trondheim, Norway}

\email{peder.thompson@ntnu.no}

\urladdr{https://folk.ntnu.no/pedertho}

\date{12 September 2019} % Activate to display a given date or no date

\thanks{L.W.C.\ was partly supported by Simons Foundation
  collaboration grant 428308.
  S.E.\ was partly supported by grant MTM2016-77445-P and FEDER funds
  and by grant 19880/GERM/15 from the Fundaci\'on S\'eneca-Agencia
  de Ciencia y Tecnolog\'{\i}a de la Regi\'on de Murcia.}



\dedicatory{To S.K.\,Jain on the occasion of his eightieth birthday}

\begin{document}

\begin{abstract}
  We introduce a notion of total acyclicity associated to a
  subcategory of an abelian category and consider the Gorenstein
  objects they define. These Gorenstein objects form a Frobenius
  category, whose induced stable category is equivalent to the
  homotopy category of totally acyclic complexes. Applied to the
  flat--cotorsion theory over a coherent ring, this provides a new
  description of the category of cotorsion Gorenstein flat modules;
  one that puts it on equal footing with the category of Gorenstein
  projective modules.
\end{abstract}

\maketitle

%%%%%%%%%%%%%%%%%%%%%%%%%%%%%%%%%%%%%%%%%%%%%%%%%%%%%%%%%%%%%%%%%%%%%%
\section*{Introduction}

\noindent
Let $A$ be an associative ring. It is classic that the stable category
of Gorenstein projective $A$-modules is triangulated equivalent to the
homotopy category of totally acyclic complexes of projective
$A$-modules. Under extra assumptions on $A$ this equivalence can be
found already in Buchweitz's 1986 manuscript \cite{ROB86}. In this
paper we focus on a corresponding equivalence for Gorenstein flat
modules. It could be pieced together from results in the literature,
but we develop a framework that provides a direct proof while also
exposing how closely the homotopical behavior of cotorsion Gorenstein
flat modules parallels that of Gorenstein projective modules.

The category of Gorenstein flat $A$-modules is rarely Frobenius, indeed we
prove in \thmref{perfect} that it only happens when every module is
cotorsion. This is evidence that one should restrict attention to the
category of cotorsion Gorenstein flat $A$-modules; in fact, it is already
known from work of Gillespie \cite{JGl17} that this category is
Frobenius if $A$ is coherent. The associated stable category is
equivalent to the homotopy category of \Ftac\ complexes of
flat-cotorsion $A$-modules; this follows from a theorem by Estrada and
Gillespie \cite{SEsJGl19} combined with recent work of Bazzoni,
Cort\'es Izurdiaga, and Estrada \cite{BCE}. The proof in
\cite{SEsJGl19} involves model structures on categories of complexes
of projective modules, and one goal of this paper---with a view
towards extending the result to non-affine schemes \cite{CET}---is to
give a proof that avoids projective modules; we achieve this with
\corref{eq0}.

The pure derived category of flat $A$-modules is the Verdier quotient
of the homotopy category of complexes of flat $A$-modules by the
subcategory of pure-acyclic complexes; its subcategory of \Ftac\
complexes was studied by Murfet and Salarian~\cite{DMfSSl11}. We show
in \thmref{I} that this subcategory is equivalent to the homotopy
category of \Ftac\ complexes of flat-cotorsion $A$-modules, and thus
to the stable category of cotorsion Gorenstein flat modules. Combining
this with results of Christensen and Kato \cite{LWCKKt18} and Estrada
and Gillespie \cite{SEsJGl19}, one can derive that under extra
assumptions on $A$, made explicit in \corref{GPGF}, the stable
category of Gorenstein projective $A$-modules is equivalent to the
stable category of cotorsion Gorenstein flat $A$-modules.

Underpinning the results we have highlighted above are a framework,
developed in Sections 1--3, and two results, \thmref[Theorems~]{2Ftac}
and \thmref[]{GG}, that show---as the semantics might suggest---that
the cotorsion Gorenstein flat modules are, indeed, the Gorenstein
modules naturally attached to the flat--cotorsion theory.
\begin{equation*}
  \ast \ \ \ast \ \ \ast
\end{equation*}
Let $\sfA$ be an abelian category and $\sfU$ a subcategory of
$\sfA$. In \dfnref[]{tac} we define a \emph{right $\sfU$-totally
  acyclic complex} to be an acyclic $\cathom{\sfA}{-}{\capUU}$-exact
complex of objects from $\sfU$ with cycle objects in $\sfUperp$. Left
$\sfU$-total acyclicity is defined dually, and in the case of a
self-orthogonal subcategory, left and right total acyclicity is the
same; see \prpref{W}. These definitions recover the standard notions
of totally acyclic complexes of projective or injective objects; see
\exaref{tac}.  In the context of a cotorsion pair $\sfUV$ the natural
complexes to consider are right $\sfU$-totally acyclic complexes, left
$\sfV$-totally acyclic complexes, and $(\capUV)$-totally acyclic
complexes for the self-orthogonal category $\capUV$.

In \secref{gor} we define left and right $\sfU$-Gorenstein objects to
be cycles in left and right $\sfU$-totally acyclic complexes. In the
context of a cotorsion pair $\sfUV$, we show that the categories of
right $\sfU$-Gorenstein objects and left $\sfV$-Gorenstein objects are
Frobenius categories whose projective-injective objects are those in
$\capUV$; see \thmref[Theorems~]{Frobenius-R} and
\thmref[]{Frobenius-L}.  In \secref{equiv} the stable categories
induced by these Frobenius categories are shown to be equivalent to
the corresponding homotopy categories of totally acyclic complexes. In
particular, \corref{s-o} recovers the classic results for Gorenstein
projective objects and Gorenstein injective objects.

The literature contains a variety of generalized notions of totally
acyclic complexes and Gorenstein objects; see for example
Sather-Wagstaff, Sharif, and White \cite{SSW-08}. We make detailed
comparisons in \rmkref{compare}; at this point it suffices to say that
our notion of Gorensteinness differs from the existing generalizations
by exhibiting periodicity: For a self-orthogonal category $\sfW$, the
category of ($\sfW$-Gorenstein)-Gorenstein objects is simply $\sfW$;
see \prpref{square}.

%%%%%%%%%%%%%%%%%%%%%%%%%%%%%%%%%%%%%%%%%%%%%%%%%%%%%%%%%%%%%%%%%%%%%%
\section{Total acyclicity and other terminology}
\label{sec:tac}

\noindent
Throughout this paper, $\sfA$ denotes an abelian category; we write
$\operatorname{hom_\sfA}$ for the hom-sets and the induced functor
from $\sfA$ to abelian groups. Tacitly, subcategories of $\sfA$ are
assumed to be full and closed under isomorphisms. A subcategory of
$\sfA$ is called \emph{additively closed} if it is additive and closed
under direct summands.

A complex of objects from $\sfA$ is referred to as an $\sfA$-complex.
We use homological notation for complexes, so for a complex $T$ the
object in degree $i$ is denoted $T_i$ and $\Cy[i]{T}$ denotes the
cycle subobject in degree $i$.

Let $\sfU$ and $\sfV$ be subcategories of $\sfA$. The \emph{right
  orthogonal} of $\sfU$ is the subcategory
\begin{equation*}
  \sfU^\perp=\setof{N\in\sfA}{\Ext[\sfA]{1}{U}{N}=0
    \text{ for all $U\in \sfU$}}\:;
\end{equation*}
the \emph{left orthogonal} of $\sfV$ is the subcategory
\begin{equation*}
  ^\perp\sfV=\setof{M\in\sfA}{\Ext[\sfA]{1}{M}{V}=0
    \text{ for all $V\in \sfV$}}\:.
\end{equation*}
In case $\sfU^\perp=\sfV$ and $^\perp\sfV=\sfU$ hold, the pair
$(\sfU,\sfV)$ is referred to as a \emph{cotorsion pair}.

In this section and the next, we develop notions of total acyclicity,
and corresponding notions of Gorenstein objects, associated to any
subcategory of $\sfA$.  Our primary applications are in the context of
a cotorsion pair.

\begin{dfn}
  \label{dfn:tac}
  Let $\sfU$ and $\sfV$ be subcategories of $\sfA$.
  \begin{prt}
  \item[\rm (\scR)] An $\sfA$-complex $T$ is called \emph{right
      $\sfU$-totally acyclic} if the following hold:
    \begin{rqm}
    \item $T$ is acyclic.
    \item For each $i\in\ZZ$ the object $T_i$ belongs to $\sfU$.
    \item For each $i\in\ZZ$ the object $\Cy[i]{T}$ belongs to
      $\sfUperp$.
    \item For each $W \in \capUU$ the complex $\cathom{\sfA}{T}{W}$ is
      acyclic.
    \end{rqm}

  \item[\rm (\scL)] An $\sfA$-complex $T$ is called \emph{left
      $\sfV$-totally acyclic} if the following hold:
    \begin{rqm}
    \item $T$ is acyclic.
    \item For each $i\in\ZZ$ the object $T_i$ belongs to $\sfV$.
    \item For each $i\in\ZZ$ the object $\Cy[i]{T}$ belongs to
      $\sfperpV$.
    \item For each $W \in \capVV$ the complex $\cathom{\sfA}{W}{T}$ is
      acyclic.
    \end{rqm}
  \end{prt}
\end{dfn}

\begin{exa}
  \label{exa:Wtac}
  Let $\sfU$ and $\sfV$ be subcategories of $\sfA$. For every
  $W \in \capUU$ a complex of the form $0 \lra W \xra{=} W \lra 0$ is
  right $\sfU$-totally acyclic; similarly, for every $W \in \capVV$
  such a complex is left $\sfV$-totally acyclic.
\end{exa}

\begin{prp}
  \label{prp:tac}
  Let $\sfU$ and $\sfV$ be subcategories of $\sfA$.
  \begin{prt}
  \item[\rm (\scR)] An $\sfA$-complex $T$ is right $\sfU$-totally
    acyclic if and only if the following hold:
    \begin{rqm}
    \item $T$ is acyclic.
    \item For each $i\in\ZZ$ the object $T_i$ belongs to $\capUU$.
    \item For each $U \in \sfU$ the complex $\cathom{\sfA}{U}{T}$ is
      acyclic.
    \item For each $W \in \capUU$ the complex $\cathom{\sfA}{T}{W}$ is
      acyclic.
    \end{rqm}

  \item[\rm (\scL)] An $\sfA$-complex $T$ is left $\sfV$-totally
    acyclic if and only if the following hold:
    \begin{rqm}
    \item $T$ is acyclic.
    \item For each $i\in\ZZ$ the object $T_i$ belongs to $\capVV$.
    \item For each $V \in \sfV$ the complex $\cathom{\sfA}{T}{V}$ is
      acyclic.
    \item For each $W \in \capVV$ the complex $\cathom{\sfA}{W}{T}$ is
      acyclic.
    \end{rqm}
  \end{prt}
\end{prp}

\begin{prf*}
  (\scR): A complex $T$ that satisfies \dfnref{tac}(\scR) trivially
  satisfies conditions $(1)$, $(3)$, and $(4)$, while $(2)$ follows
  from \dfnref[]{tac}(\scR.2) and \dfnref[]{tac}(\scR.3) as $\sfUperp$
  is closed under extensions. Conversely, a complex $T$ that satisfies
  conditions $(1)$--$(4)$ in the statement trivially satisfies
  conditions $(1)$, $(2)$, and $(4)$ in \dfnref{tac}(\scR).  Moreover
  it follows from (2) and (3) that also condition
  \dfnref[]{tac}(\scR.3) is satisfied.

  The proof of (\scL) is similar.
\end{prf*}

\begin{exa}
  \label{exa:tac-0}
  A right $\sfA$-totally acyclic complex is a contractible complex of
  injective objects, and a left $\sfA$-totally acyclic complex is a
  contractible complex of projective objects.
\end{exa}

In this paper we call a subcategory $\sfW$ of $\sfA$
\emph{self-orthogonal} if $\Ext[\sfA]{1}{W}{W'}=0$ holds for all $W$
and $W'$ in $\sfW$.

\begin{prp}
  \label{prp:W}
  Let $\sfW$ be a subcategory of $\sfA$. The following conditions are
  equivalent
  \begin{eqc}
  \item $\sfW$ is self-orthogonal.
  \item Every object in $\sfW$ belongs to $\sfW^\perp$.
  \item Every object in $\sfW$ belongs to ${}^\perp\sfW$.
  \item One has $\sfW \cap \sfW^\perp = \sfW = {}^\perp\sfW\cap\sfW$.
  \end{eqc}
  Moreover, if $\sfW$ satisfies these conditions, then an
  $\sfA$-complex is right $\sfW$-totally acyclic if and only if it is
  left $\sfW$-totally acyclic.
\end{prp}

\begin{prf*}
  Evidently, \eqclbl{i} implies \eqclbl{iv}, and \eqclbl{iv} implies
  both \eqclbl{ii} and \eqclbl{iii}.  Conditions \eqclbl{ii} and
  \eqclbl{iii} each precisely mean that $\Ext[\sfA]{1}{W}{W'}=0$ holds
  for all $W$ and $W'$ in $\sfW$, so either implies \eqclbl{i}.

  Now assume that $\sfW$ satisfies \eqclbl{i}--\eqclbl{iv}. Parts (1)
  are the same in \prpref{tac}(\scR) and \prpref[]{tac}(\scL), and so
  are parts (2) per the assumption
  $\sfW \cap \sfW^\perp = {}^\perp\sfW\cap\sfW$. Part (3) in
  \prpref[]{tac}(\scR) coincides with part (4) in \prpref[]{tac}(\scL)
  by the assumption $\sfW = {}^\perp\sfW\cap\sfW$, and
  \prpref[]{tac}(\scR.4) coincides with \prpref[]{tac}(\scL.3) per the
  assumption $\sfW \cap \sfW^\perp = \sfW$.
\end{prf*}

\begin{dfn}
  \label{dfn:s-o}
  For a self-orthogonal subcategory $\sfW$ of $\sfA$, a right,
  equivalently left, $\sfW$-totally acyclic complex is simply called a
  $\sfW$-\emph{totally acyclic complex.}
\end{dfn}

\begin{exa}
  \label{exa:tac}
  The subcategory $\catPrj{\sfA}$ of projective objects in $\sfA$ is
  self-orthogonal, and a $\catPrj{\sfA}$-totally acyclic complex is
  called a \emph{totally acyclic complex of projective objects}. In
  the special case where $\sfA$ is the category $\catMod{A}$ of
  modules over a ring $A$ these were introduced by Auslander and
  Bridger~\cite{MAsMBr69}; see also Enochs and
  Jenda~\cite{EEnOJn95b}. The terminology is due to Avramov and
  Martsinkovsky \cite{LLAAMr02}.

  Dually, $\catInj{\sfA}$ is the subcategory of injective objects in
  $\sfA$, and an $\catInj{\sfA}$-totally acyclic complex is called a
  \emph{totally acyclic complex of injective objects;} see Krause
  \cite{HKr05}. The case $\sfA = \catMod{A}$ was first considered
  in~\cite{EEnOJn95b}.
\end{exa}

\begin{rmk}
  \label{rmk:UVW}
  For a cotorsion pair $(\sfU,\sfV)$ in $\sfA$, the subcategory
  $\capUV$ is self-orthogonal.  It follows from \prpref{tac} that
  every right $\sfU$-totally acyclic complex and every left
  $\sfV$-totally acyclic complex is $(\capUV)$-totally acyclic.
\end{rmk}

%%%%%%%%%%%%%%%%%%%%%%%%%%%%%%%%%%%%%%%%%%%%%%%%%%%%%%%%%%%%%%%%%%%%%%
\section{Gorenstein objects}\label{sec:gor}

\noindent
In line with standard terminology, cycles in totally acyclic complexes
are called Gorenstein objects.

\begin{dfn}
  \label{dfn:Gor}
  Let $\sfU$ and $\sfV$ be subcategories of $\sfA$.
  \begin{prt}
  \item[(\scR)] An object $M$ in $\sfA$ is called \emph{right
      $\sfU$-Gorenstein} if there is a right $\sfU$-totally acyclic
    complex $T$ with $\Cy[0]{T} = M$. Denote by $\catGR{\sfA}$ the
    full subcategory of right $\sfU$-Gorenstein objects in $\sfA$.

  \item[(\scL)] An object $M$ in $\sfA$ is called \emph{left
      $\sfV$-Gorenstein} if there is a left $\sfV$-totally acyclic
    complex $T$ with $\Cy[0]{T} = M$. Denote by $\catGL{\sfA}$ the
    full subcategory of left $\sfV$-Gorenstein objects in $\sfA$.
  \end{prt}

  \noindent For a self-orthogonal subcategory $\sfW$ one has
  $\catGR[\sfW]{\sfA} = \catGL[\sfW]{\sfA}$, see \prpref{W}; this
  category is denoted $\catG{\sfA}$, and its objects are called
  $\sfW$-\emph{Gorenstein.}
\end{dfn}

Notice that if $\sfU$ is an additive subcategory, then so is
$\catGR{\sfA}$; similarly for $\sfV$ and $\catGL{\sfA}$.

\begin{exa}
  \label{exa:WGor}
  Let $\sfU$ and $\sfV$ be subcategories of $\sfA$. Objects in
  $\capUU$ are right $\sfU$-Gorenstein and objects in $\capVV$ are
  left $\sfV$-Gorenstein; see \exaref{Wtac}.
\end{exa}

\begin{rmk}
  \label{rmk:compare}
  We compare our definitions of total acyclicity and Gorenstein
  objects with others that already appear in the literature.
  \begin{enumerate}
  \item For an additive category $\sfW$, Iyengar and Krause
    \cite{SInHKr06} define a ``totally acyclic complex over $\sfW$.''
    For additive subcategories $\sfU$ and $\sfV$ of an abelian
    category, a right $\sfU$-totally acyclic complex is totally
    acyclic over $\capUU$ in the sense of \dfncite[5.2]{SInHKr06}, and
    a left $\sfV$-totally acyclic complex is totally acyclic over
    $\capVV$. In particular, for a self-orthogonal additive
    subcategory $\sfW$ of an abelian category, a $\sfW$-totally
    acyclic complex is the same as an acyclic complex that is totally
    acyclic over $\sfW$ in the sense of
    \dfncite[5.2]{SInHKr06}.
  \item For an additive subcategory $\sfW$ of an abelian category,
    Sather-Wagstaff, Sharif, and White \cite{SSW-08} define a
    ``totally $\sfW$-acyclic'' complex. A right or left $\sfW$-totally
    acyclic complex is totally $\sfW$-acyclic in the sense of
    \dfncite[4.1]{SSW-08}; the converse holds if $\sfW$ is
    self-orthogonal. For a self-orthogonal additively closed
    subcategory $\sfW$ of a module category, Geng and Ding
    \cite{YGnNDn11} study the associated Gorenstein objects.
  \item For subcategories $\sfU$ and $\sfV$ of $\catMod{A}$ with
    $\catPrj{A}\subseteq \sfU$ and $\catInj{A}\subseteq \sfV$, Pan and
    Cai \cite{PQnCFq15} define ``$(\sfU,\sfV)$-Gorenstein
    projective/injective'' modules.  In this setting, a right
    $\sfU$-Gorenstein module is $(\sfU,\capUU)$-Gorenstein projective
    in the sense of \dfncite[2.1]{PQnCFq15}, and a left
    $\sfV$-Gorenstein module is $(\capVV,\sfV)$-Gorenstein injective
    in the sense of \dfncite[2.2]{PQnCFq15}.
  \item For a complete hereditary cotorsion pair $(\sfU,\sfV)$ in an
    abelian category, Yang and Chen \cite{XYnWCh17} define a
    ``complete $\sfU$-resolution.'' Every right $\sfU$-totally acyclic
    complex is a complete $\sfU$-resolution in the sense of
    \dfncite[3.1]{XYnWCh17}.
  \item For a pair of subcategories $(\sfU,\sfV)$ in an abelian
    category, Becerril, Mendoza, and Santiago \cite{BMS} define a
    ``left complete $(\sfU,\sfV)$-resolution.'' If $(\sfU,\sfV)$ is a
    cotorsion pair, then a right $\sfU$-totally acyclic complex is a
    left complete $(\sfU,\sfU\cap\sfV)$-resolution in the sense of
    \dfncite[3.2]{BMS}.
  \end{enumerate}
  The key difference between \dfnref{tac} and those cited above is
  that \dfnref[]{tac}---motivated by \ref{cot-cyc}---places
  restrictions on the cycle objects in a totally acyclic complex; the
  significance of this becomes apparent in \prpref{square}.
\end{rmk}

\begin{rmk}
  \label{rmk:R=LGor}
  Given a cotorsion pair $\sfUV$ in $\sfA$, it follows from
  \rmkref{UVW} that there are containments
  \begin{equation*}
    \catGR{\sfA} \subseteq \catG[\sfU\cap\sfV]{\sfA}
    \supseteq \catGL{\sfA}\:.
  \end{equation*}
\end{rmk}

\begin{exa}
  \label{exa:Gor}
  A right $\sfA$-Gorenstein object is injective, and a left
  $\sfA$-Gorenstein object is projective; see \exaref{tac-0}.

  The subcategory $\catPrj{\sfA}$ is self-orthogonal, and a
  $\catPrj{\sfA}$-Gorenstein object is called \emph{Gorenstein
    projective}; see \cite{MAsMBr69,EEnOJn95b} for the special case
  $\sfA = \catMod{A}$. Similarly, an $\catInj{\sfA}$-Gorenstein object
  is called \emph{Gorenstein injective;} see \cite{HKr05} and see
  \cite{EEnOJn95b} for the case $\sfA = \catMod{A}$.
\end{exa}

The next three results, especially \prpref{square}, are motivated in
part by \cite[Theorem A]{SSW-08}. We consider what happens when one
iterates the process of constructing Gorenstein objects. Starting from
a self-orthogonal additively closed subcategory, our construction
iterated twice returns the original subcategory. The construction in
\cite{SSW-08} is, in contrast, ``idempotent.''

\begin{lem}
  \label{lem:cap}
  Let $\sfU$ and $\sfV$ be additively closed subcategories of
  $\sfA$. One has
  \begin{align*}
    {}^\perp\catGR{\sfA} \cap \catGR{\sfA}
    &\deq \capUU
      \deq \catGR{\sfA} \cap \catGR{\sfA}^\perp\quad\text{and}\\
    {}^\perp\catGL{\sfA} \cap \catGL{\sfA}
    &\deq \capVV
      \deq \catGL{\sfA} \cap \catGL{\sfA}^\perp\:.
  \end{align*}
  In particular, $\catGR{\sfA}$ is self-orthogonal if and only if
  $\catGR{\sfA} = \capUU$ holds, and $\catGL{\sfA}$ is self-orthogonal
  if and only if $\catGL{\sfA} = \capVV$ holds.

  For a self-orthogonal category $\sfW$ one has
  \begin{equation}
    \label{eq:w}
    {}^\perp\catG{\sfA} \cap \catG{\sfA} = \sfW
    = \catG{\sfA} \cap \catG{\sfA}^\perp\:.
  \end{equation}
\end{lem}

\begin{prf*}
  Set $\sfW = \capUU$ and notice that $\sfW$ is self-orthogonal and
  additively closed.  By \exaref{WGor} objects in $\sfW$ are right
  $\sfU$-Gorenstein, and by \prpref{tac} the subcategory $\sfW$ is
  contained in both ${}^\perp\catGR{\sfA}$ and
  $\catGR{\sfA}^\perp$. Let $G$ be a right $\sfU$-Gorenstein
  object. By \prpref{tac} there are exact sequences
  \begin{equation*}
    \eta' \deq 0 \to G' \to T' \to G \to 0
    \qqand \eta'' \deq 0 \to G \to T'' \to G'' \to 0
  \end{equation*}
  where $G'$ and $G''$ are right $\sfU$-Gorenstein, while $T'$ and
  $T''$ belong to $\sfW$. If $G$ belongs to ${}^\perp\catGR{\sfA}$,
  then $\eta'$ splits, so $G$ is a summand of $T'$ and hence in
  $\sfW$. Similarly, if $G$ is in $\catGR{\sfA}^\perp$, then $\eta''$
  splits, and it follows that $G$ is in $\sfW$. This proves the first
  set of equalities, and the ones pertaining to $\catGL{\sfA}$ are
  proved similarly.

  The remaining assertions are immediate in view of \prpref{W}.
\end{prf*}

\begin{rmk}
  Let $\sfU$ be an additively closed subcategory of $\sfA$. It follows
  from \exaref{WGor} and \lemref{cap} that objects in $\capUU$ are
  both right $\sfU$-Gorenstein and right $\catGR{\sfA}$-Gorenstein. On
  the other hand, a right $\catGR{\sfA}$-Gorenstein object belongs by
  \dfnref{tac}(\scR.3) to $\catGR{\sfA}^\perp$, so any object that is
  both right $\sfU$- and right $\catGR{\sfA}$-Gorenstein belongs to
  $\capUU$. In symbols,
  \begin{equation*}
    \catGR{\sfA} \cap \mathsf{RGor}_{\catGR{\sfA}}(\sfA) \deq \capUU\:.
  \end{equation*}
  For an additively closed subcategory $\sfV$, a similar argument
  yields
  \begin{equation*}
    \catGL{\sfA} \cap \mathsf{LGor}_{\catGL{\sfA}}(\sfA) \deq \capVV\:.
  \end{equation*}
\end{rmk}

\begin{prp}
  \label{prp:square}
  Let $\sfW$ be a self-orthogonal additively closed subcategory of
  $\sfA$.  A right or left $\catG{\sfA}$-totally acyclic complex is a
  contractible complex of objects from $\sfW$. In particular, one has
  \begin{equation}
    \label{eq:square-1}
    \catGL[\mathsf{Gor}_{\sfW}(\sfA)]{\sfA} = \sfW =
    \catGR[\mathsf{Gor}_{\sfW}(\sfA)]{\sfA}\:.
  \end{equation}
  Moreover, the following hold
  \begin{itemlist}
  \item If\, $(\catG{\sfA},\catG{\sfA}^\perp)$ is a cotorsion pair,
    then one has
    \begin{equation}
      \label{eq:square-3}
      \catGL[\mathsf{Gor}_{\sfW}(\sfA)^\perp]{\sfA} =
      \catG{\sfA}\:.
    \end{equation}
  \item If\, $({}^\perp\catG{\sfA},\catG{\sfA})$ is a cotorsion pair,
    then one has
    \begin{equation}
      \label{eq:square-4}
      \catGR[{}^\perp\mathsf{Gor}_{\sfW}(\sfA)]{\sfA} = \catG{\sfA}\:.
    \end{equation}
  \end{itemlist}
\end{prp}

\begin{prf*}
  A right $\catG{\sfA}$-totally acyclic complex $T$ is by \prpref{tac}
  and \eqref{w} an acyclic complex of objects from $\sfW$, and by
  \dfnref{tac} the cycles $\Cy[i]{T}$ belong to
  $\catG{\sfA}^\perp$. As $\sfW$ is contained in $\catG{\sfA}$, it
  follows from \prpref{tac}(\scR.3) that the cycles $\Cy[i]{T}$ are
  contained in $\sfW^\perp$. It now follows from \dfnref{tac} that $T$
  is $\sfW$-totally acyclic, whence the cycles $\Cy[i]{T}$ belong to
  $\catG{\sfA}$ and hence to $\sfW$, see \eqref{w}. Thus $T$ is an
  acyclic complex of objects from $\sfW$ with cycles in
  $\sfW \subset \sfW^\perp$ and, therefore, contractible. A parallel
  argument shows that a left $\catG{\sfA}$-totally acyclic complex is
  contractible.

  Assume that $(\catG{\sfA},\catG{\sfA}^\perp)$ is a cotorsion pair;
  by \eqref{w} and \rmkref{R=LGor} one has
  $\catGL[\mathsf{Gor}_{\sfW}(\sfA)^\perp]{\sfA} \subseteq
  \catG{\sfA}$. To prove the opposite containment, let $T$ be a
  $\sfW$-totally acyclic complex. By \dfnref{tac} it is an acyclic
  complex of objects from $\sfW \subseteq \catG{\sfA}^\perp$, see
  \eqref{w}, and $\cathom{\sfA}{W}{T}$ is acyclic for every object $W$
  in $\catG{\sfA} \cap \catG{\sfA}^\perp$. Moreover, the cycle objects
  $\Cy[i]{T}$ belong to $\catG{\sfA}$ by \dfnref{Gor}, so $T$ is per
  \dfnref{tac} a left $\catG{\sfA}^\perp$-totally acyclic complex.

  A parallel argument proves the last assertion.
\end{prf*}

\begin{exa}
  Let $A$ be a ring. \v{S}aroch and
  {\v{S}}tov{\'{\i}}{\v{c}}ek~\thmcite[4.6]{JSrJSt} show that the
  subcategory $\catG[\mathsf{Inj}]{A}$ of Gorenstein injective
  $A$-modules is the right half of a cotorsion pair, so by
  \prpref{square} one has
  $\catGL[\mathsf{Gor_{Inj}}(A)]{A} = \catInj{A} =
  \catGR[\mathsf{Gor_{Inj}}(A)]{A}$ and
  $\catGR[{}^\perp\mathsf{Gor_{Inj}}(A)]{A} = \catG[\mathsf{Inj}]{A}$.
\end{exa}

\begin{lem}
  \label{lem:extclosed}
  Let $\sfU$ and $\sfV$ be additive subcategories of $\sfA$.
  \begin{prt}
  \item[\rm (\scR)] The subcategory \smash{$\catGR{\sfA}$} is closed
    under extensions.
  \item [\rm (\scL)] The subcategory $\catGL{\sfA}$ is closed under
    extensions.
  \end{prt}
\end{lem}

\begin{prf*}
  Let $0\to M' \to M \to M'' \to 0$ be an exact sequence where $M'$
  and $M''$ are right $\sfU$-Gorenstein objects. Let $T'$ and $T''$ be
  right $\sfU$-totally acyclic complexes with $\Cy[0]{T'}=M'$ and
  $\Cy[0]{T''} = M''$. Per \rmkref{compare}(1) it follows from
  \prpcite[4.4]{SSW-08} that there exists an $\sfA$-complex $T$ that
  satisfies conditions (1), (2), and (4) in \dfnref{tac}(\scR), has
  $\Cy[0]{T}=M$, and fits in an exact sequence
  \begin{equation*}
    0 \lra T' \lra T \lra T'' \lra 0\:.
  \end{equation*}
  The functor $\Cy{-}$ is left exact, and since $T'$ is acyclic a
  standard application of the Snake Lemma yields an exact sequence
  \begin{equation*}
    0 \lra \Cy[i]{T'} \lra \Cy[i]{T} \lra \Cy[i]{T''} \lra 0
  \end{equation*}
  for every $i\in\ZZ$.  As $\sfUperp$ is closed under extensions, it
  follows that $\Cy[i]{T}$ belongs to $\sfUperp$ for each $i$ and thus
  $T$ is right $\sfU$-totally acyclic by \dfnref{tac}. This proves
  (\scR) and a similar argument proves (\scL).
\end{prf*}

\begin{thm}
  \label{thm:Frobenius-R}
  Let $\sfU$ be an additively closed subcategory of $\sfA$.  The
  category $\catGR{\sfA}$ is Frobenius and $\capUU$ is the subcategory
  of projective-injective objects.
\end{thm}

\begin{prf*}
  Set $\sfW = \capUU$ and notice that $\sfW$ is additively closed. It
  follows from \lemref{extclosed} that $\catGR{\sfA}$ is an exact
  category. It is immediate from \exaref{WGor} and \prpref{tac} that
  objects in $\sfW$ are both projective and injective in
  $\catGR{\sfA}$. It is now immediate from \dfnref{Gor} that
  \smash{$\catGR{\sfA}$} has enough projectives and injectives. It
  remains to show that every projective and every injective object in
  $\catGR{\sfA}$ belongs to $\sfW$.

  Let $P$ be a projective object in \smash{$\catGR{\sfA}$}. By
  \dfnref{Gor} and \prpref{tac} there is an exact sequence
  $0 \to P' \to W \to P \to 0$ in $\sfA$ with $P'\in\catGR{\sfA}$ and
  $W \in \sfW$. As all three objects belong to \smash{$\catGR{\sfA}$}
  it follows by projectivity of $P$ that the sequence splits, so $P$
  is a summand of $W$, and thus in $\sfW$. A dual argument shows that
  every injective object in $\catGR{\sfA}$ belongs to $\sfW$. Thus
  $\catGR{\sfA}$ is a Frobenius category and $\sfW$ is the subcategory
  of projective-injective objects.
\end{prf*}

\begin{thm}
  \label{thm:Frobenius-L}
  Let $\sfV$ be an additively closed subcategory of $\sfA$.  The
  category $\catGL{\sfA}$ is Frobenius and $\capVV$ is the subcategory
  of projective-injective objects.
\end{thm}

\begin{prf*}
  Parallel to the proof of \thmref{Frobenius-R}.
\end{prf*}

%%%%%%%%%%%%%%%%%%%%%%%%%%%%%%%%%%%%%%%%%%%%%%%%%%%%%%%%%%%%%%%%%%%%%%
\section{An equivalence of triangulated categories}
\label{sec:equiv}

\noindent
Generalizing the classic result, we prove here that the stable
category of right/left Gorenstein objects is equivalent to the
homotopy category of right/left totally acyclic complexes.

\begin{lem}
  \label{lem:lift0}
  Let $\sfU$ be a subcategory of $\sfA$; let $T$ and $T'$ be right
  $\sfU$-totally acyclic complexes. Every morphism
  $\mapdef{\grf}{\Cy[0]{T}}{\Cy[0]{T'}}$ in $\sfA$ lifts to a morphism
  $\mapdef{\phi}{T}{T'}$ of $\sfA$-complexes.
\end{lem}

\begin{prf*}
  Let a morphism $\mapdef{\grf}{\Cy[0]{T}}{\Cy[0]{T'}}$ be given; to
  see that it lifts to a morphism $\mapdef{\phi}{T}{T'}$ of complexes
  it is sufficient to show that $\grf$ lifts to morphisms
  $\mapdef{\phi_1}{T_1}{T'_1}$ and $\mapdef{\phi_0}{T_0}{T'_0}$. As
  $T_1$ is in $\sfU$ and $T'$ is right $\sfU$-totally acyclic, one
  obtains per \prpref{tac}(\scR.3) an exact sequence
  \begin{equation*}
    0 \lra \cathom{\sfA}{T_1}{\Cy[1]{T'}}
    \lra \cathom{\sfA}{T_1}{T'_1}
    \lra \cathom{\sfA}{T_1}{\Cy[0]{T'}} \lra 0\:.
  \end{equation*}
  In particular, there is a $\phi_1\in\cathom{\sfA}{T_1}{T'_1}$ with
  $\partial^{T'}_1\phi_1 = \grf\partial_1^T$. As $T'_0$ is in $\capUU$
  and $T$ is right $\sfU$-totally acyclic, it follows that the
  sequence
  \begin{equation*}
    0 \lra \cathom{\sfA}{\Cy[-1]{T}}{T'_0}
    \lra \cathom{\sfA}{T_0}{T'_0} \lra \cathom{\sfA}{\Cy[0]{T}}{T'_0} \lra 0
  \end{equation*}
  is exact, whence there exists a $\phi_0\in\cathom{\sfA}{T_0}{T'_0}$
  that lifts $\grf$.
\end{prf*}

\begin{lem}
  \label{lem:lift}
  Let $\sfU$ be a subcategory of $\sfA$ and $\mapdef{\phi}{T}{T'}$ be
  a morphism of right $\sfU$-totally acyclic complexes. If the cycle
  subobject $\Cy[0]{T}$ has a decomposition
  $\Cy[0]{T} = Z \oplus \widetilde{Z}$ with $Z \subseteq \ker{\phi_0}$
  and $\widetilde{Z}\in\sfU$, then $\phi$ is null-homotopic.
\end{lem}

\begin{prf*}
  The goal is to construct a family of morphisms
  $\mapdef{\grs_i}{T_i}{T'_{i+1}}$ such that
  $\phi_i = \partial_{i+1}^{T'}\grs_i + \grs_{i-1}\partial_i^T$ holds
  for all $i\in\ZZ$. Set $\widetilde\grf =
  \phi_0|_{\widetilde{Z}}$. By \dfnref{tac} each object $\Cy[i]{T'}$
  is in $\sfU^\perp$. It follows that there is an exact sequence,
  \begin{equation*}
    0 \lra \cathom{\sfA}{\widetilde{Z}}{\Cy[1]{T'}}
    \lra \cathom{\sfA}{\widetilde{Z}}{T'_1}
    \lra \cathom{\sfA}{\widetilde{Z}}{\Cy[0]{T'}}
    \lra 0\:.
  \end{equation*}
  In particular, there is a
  $\widetilde\grs \in \cathom{\sfA}{\widetilde{Z}}{T'_1}$ with
  $\partial_1^{T'}\widetilde\grs = \widetilde\grf$. Set
  $\widetilde\grs_0 = 0\oplus\widetilde\grs$; by exactness of the
  sequence
  \begin{equation*}
    0 \lra \cathom{\sfA}{\Cy[-1]{T}}{T'_1}
    \lra \cathom{\sfA}{T_0}{T'_1}
    \lra \cathom{\sfA}{\Cy[0]{T}}{T'_1}
    \lra 0\:,
  \end{equation*}
  $\widetilde\grs_0$ lifts to a morphism $\mapdef{\grs_0}{T_0}{T'_1}$.

  We proceed by induction to construct the morphisms $\grs_i$ for
  $i \ge 1$. The image of the morphism $\phi_1 - \grs_0\partial_1^T$
  is in $\Cy[1]{T'}$ as one has
  \begin{align*}
    \partial_1^{T'}(\phi_1 - \grs_0\partial_1^T) %%
    &=\phi_0\partial_1^{T} - \partial_1^{T'}\grs_0\partial_1^T \\
    &= (0 \oplus \widetilde\grf)\partial_1^{T} - \partial_1^{T'}(0
      \oplus \widetilde\grs)\partial_1^T\\
    &= (0 \oplus (\widetilde\grf
      - \partial_1^{T'}\widetilde\grs))\partial_1^T\\
    &= 0\:.
  \end{align*}
  As $T_1$ is in $\sfU$ and $\Cy[2]{T'}$ is in $\sfU^\perp$ per
  \dfnref{tac}, there is an exact sequence
  \begin{equation*}
    0 \lra \cathom{\sfA}{T_1}{\Cy[2]{T'}}
    \lra \cathom{\sfA}{T_1}{T'_2} \lra
    \cathom{\sfA}{T_1}{\Cy[1]{T'}} \lra 0\:.
  \end{equation*}
  In particular, there is a $\grs_1 \in \cathom{\sfA}{T_1}{T'_2}$ with
  $\partial_2^{T'}\grs_1 = \phi_1 - \grs_0\partial_1^T$. Now let
  $i\ge 1$ and assume that $\grs_j$ has been constructed for
  $0 \le j \le i$. The standard computation
  \begin{align*}
    \partial_{i+1}^{T'}(\phi_{i+1} - \grs_{i}\partial_{i+1}^T) %%
    & = (\phi_{i} - \partial_{i+1}^{T'}\grs_{i})\partial_{i+1}^T \\
    & = (\grs_{i-1}\partial_{i}^T)\partial_{i+1}^T \\
    & = 0\:
  \end{align*}
  shows that the image of $\phi_{i+1} - \grs_{i}\partial_{i+1}^T$ is
  in $\Cy[i+1]{T'}$. As $T_{i+1}$ is in $\sfU$ and $\Cy[i+2]{T'}$ is
  in $\sfU^\perp$, the existence of the desired $\grs_{i+1}$ follows
  as for $i=0$.

  Finally, we prove the existence of the morphisms $\grs_i$ for
  $i \le -1$ by descending induction. The morphism
  $\mapdef{\phi_0 - \partial_1^{T'}\grs_0}{T_0}{T'_0}$ restricts to
  $0$ on $\Cy[0]{T}$; indeed one has
  \begin{equation*}
    (\phi_0 - \partial_1^{T'}\grs_0)|_{\Cy[0]{T}} =
    0\oplus\widetilde\grf - \partial_1^{T'}(0\oplus \widetilde\grs) =
    0\oplus(\widetilde\grf - \partial_1^{T'}\widetilde\grs) = 0\:.
  \end{equation*}
  Thus it induces a morphism $\grz_{-1}$ from
  $T_0/\Cy[0]{T} \is \Cy[-1]{T}$ to $T'_0$ with
  $\grz_{-1}\partial_0^T = \phi_0 - \partial_1^{T'}\grs_0$.  As $T_0'$
  is in $\capUU$ it follows that the sequence
  \begin{equation*}
    0 \lra \cathom{\sfA}{\Cy[-2]{T}}{T'_0}
    \lra \cathom{\sfA}{T_{-1}}{T'_0}
    \lra \cathom{\sfA}{\Cy[-1]{T}}{T'_0} \lra 0
  \end{equation*}
  is exact. In particular, there is a
  $\grs_{-1}\in\cathom{\sfA}{T_{-1}}{T'_0}$ with
  $\grs_{-1}|_{\Cy[-1]{T}}=\grz_{-1}$ and, therefore,
  $\grs_{-1}\partial_0^T = \phi_0 - \partial_1^{T'}\grs_0$. Now let
  $i \le -1$ and assume that $\grs_j$ has been constructed for
  $0 \ge j \ge i$. The standard computation
  \begin{equation*}
    (\phi_i - \partial^{T'}_{i+1}\grs_i)\partial_{i+1}^T =
    \partial^{T'}_{i+1}(\phi_{i+1} - \grs_i\partial_{i+1}^T) =
    \partial^{T'}_{i+1}(\partial^{T'}_{i+2} \grs_{i+1}) = 0
  \end{equation*}
  shows that the morphism $\phi_i - \partial^{T'}_{i+1}\grs_i$
  restricts to $0$ on $\Cy[i]{T}$. It follows that it induces a
  morphism $\grz_{i-1}$ on $T_i/\Cy[i]{T} \is \Cy[i-1]{T}$ with
  $\grz_{i-1}\partial_i^T = \phi_i - \partial^{T'}_{i+1}\grs_i$. Since
  $T_i'$ is in $\capUU$, it follows as for $i=0$ that the desired
  $\grs_{i-1}$ exists.
\end{prf*}

\begin{prp}
  \label{prp:lift}
  Let $\sfU$ be a subcategory of $\sfA$. Let $T$ and $T'$ be right
  $\sfU$-totally acyclic complexes and
  $\mapdef{\grf}{\Cy[0]{T}}{\Cy[0]{T'}}$ be a morphism in $\sfA$.
  \begin{prt}
  \item If $\mapdef{\phi}{T}{T'}$ and $\mapdef{\psi}{T}{T'}$ are
    morphisms that lift $\grf$, then $\phi-\psi$ is null-homotopic.
  \item If $\grf$ is an isomorphism and $\mapdef{\phi}{T}{T'}$ is a
    morphism that lifts $\grf$, then $\phi$ is a homotopy equivalence.
  \end{prt}
\end{prp}

\begin{prf*}
  (a): Immediate from \lemref{lift} as
  $(\phi-\psi)|_{\Cy[0]{T}} = \grf-\grf=0$.

  (b): Let $\mapdef{\phi'}{T'}{T}$ be a lift of $\grf^{-1}$; see
  \lemref{lift0}. The restriction of $1^T - \phi'\phi$ to $\Cy[0]{T}$
  is $0$, so it follows from part (a) that $1^T - \phi'\phi$ is
  null-homotopic. Similarly, $1^{T'} - \phi\phi'$ is null-homotopic;
  that is, $\phi$ is a homotopy equivalence.
\end{prf*}

\begin{dfn}
  \label{dfn:Ktac}
  Let $\sfU$ and $\sfV$ be subcategories of $\sfA$. Denote by
  $\catKRtac{\capUU}$ and $\catKLtac{\capVV}$ the homotopy categories
  of right $\sfU$-totally acyclic complexes and left $\sfV$-totally
  acyclic complexes.

  The subcategory $\capUU$ is self-orthogonal, so the categories
  $\catKRtac[\smash{(\capUU)}]{\capUU}$ and
  $\catKLtac[\smash{(\capUU)}]{\capUU}$ coincide, see \prpref{W}, and
  are denoted $\catKtac{\capUU}$.  The self-orthogonal subcategory
  $\capVV$ similarly gives a category $\catKtac{\capVV}$.  For a
  cotorsion pair $\sfUV$ all of these homotopy categories are
  $\catKtac{\capUV}$.
\end{dfn}

If $\sfU$ is an additive subcategory, then the homotopy category
$\catKRtac{\capUU}$ is triangulated; similarly for $\sfV$ and
$\catKLtac{\capVV}$.

\begin{lem}
  \label{lem:split}
  Let $\sfU$ be an additively closed subcategory of $\sfA$.  Let $T$
  be a right $\sfU$-totally acyclic complex; if $\Cy[i]{T}$ belongs to
  $\sfU$ for some $i\in\ZZ$, then $T$ is contractible.
\end{lem}

\begin{prf*}
  Set $\sfW = \capUU$ and notice that $\sfW$ is additively closed.  To
  prove that $T$ is contractible it is enough to show that
  $Z_i := \Cy[i]{T}$ belongs to $\sfW$ for every $i\in\ZZ$. There are
  exact sequences
  \begin{equation*}
    \tag{$\ast$}
    0 \lra  Z_{j+1} \lra T_{j+1} \lra Z_j \lra 0
  \end{equation*}
  with $T_{j+1}$ in $\sfW$ and $Z_{j+1},Z_j \in \sfU^\perp$; see
  \dfnref{tac} and \prpref{tac}. Without loss of generality, assume
  that $Z_0$ is in $\sfU$ and hence in $\sfW$.

  Let $j\ge 0$ and assume that $Z_j$ is $\sfW$. The sequence $(\ast)$
  splits as $Z_{j}$ is in $\sfU$ and $Z_{j+1}$ is in $\sfU^\perp$. It
  follows that $Z_{j+1}$ is in $\sfW$, so by induction $Z_i$ is in
  $\sfW$ for all $i\ge 0$.

  Now let $j<0$ and assume that $Z_{j+1}$ is in $\sfW$. The sequence
  $(\ast)$ splits as $\cathom{\sfA}{T}{Z_{j+1}}$ is acyclic by
  \dfnref{tac}. It follows that $Z_j$ belongs to $\sfW$, so by
  descending induction $Z_i$ is in $\sfW$ for all $i\le 0$.
\end{prf*}

\begin{prp}
  \label{prp:Tdot}
  Let $\sfU$ be an additively closed subcategory of $\sfA$.
  \begin{itemlist}
  \item For every right $\sfU$-Gorenstein object $M$ fix a right
    $\sfU$-totally acyclic complex $T$ with $\Cy[0]{T} = M$ and set
    $\fuTdot(M) = T$.
  \item For every morphism $\mapdef{\grf}{M}{M'}$ of right
    $\sfU$-Gorenstein objects fix by {\rm \lemref[]{lift0}} a lift
    $\mapdef{\phi}{\fuTdot(M)}{\fuTdot(M')}$ of $\grf$ and set
    $\fuTdot(\grf) = [\phi]$.
  \end{itemlist}
  This defines a functor
  \begin{equation*}
    \mapdef[\lra]{\fuTdot}{\catGR{\sfA}}{\catKRtac{\capUU}}\:.
  \end{equation*}
  For every morphism $\grf$ in $\catGR{\sfA}$ that factors through an
  object in $\capUU$ one has $\fuTdot(\grf) = [0]$. In particular,
  $\fuTdot(M)$ is contractible for every $M$ in $\capUU$.
\end{prp}

\begin{prf*}
  Let $M$ be a right $\sfU$-Gorenstein object. Denote by
  $\mapdef{\gri^M}{\fuTdot(M)}{\fuTdot(M)}$ the fixed lift of $1^M$;
  that is, $[\gri^M] = \fuTdot(1^M)$. As the morphisms
  $1^{\fuTdot(M)}$ and $\gri^M$ agree on $\Cy[0]{\fuTdot(M)}=M$, it
  follows from \lemref{lift} that the difference
  $1^{\fuTdot(M)} - \gri^M$ is null-homotopic. That is, one has
  $\fuTdot(1^M) = [1^{\fuTdot(M)}]$, which is the identity on
  $\fuTdot(M)$ in $\catKRtac{\capUU}$.

  Let $M' \xra{\grf'} M \xra{\grf} M''$ be morphisms of right
  $\sfU$-Gorenstein objects. The restrictions of $\fuTdot(\grf\grf')$
  and $\fuTdot(\grf)\fuTdot(\grf')$ to $\Cy[0]{\fuTdot(M')}$ are both
  $\grf\grf'$. It now follows from Lemma~\ref{lem:lift} that the
  homotopy classes $\fuTdot(\grf\grf')$ and
  $\fuTdot(\grf)\fuTdot(\grf')$ are equal. Thus $\fuTdot$ is a
  functor.

  For an object $M$ in the additively closed subcategory $\capUU$ it
  follows from \lemref{split} that $\fuTdot(M)$ is contractible.
  Finally, if a morphism $\mapdef{\grf}{M'}{M''}$ in $\catGR{\sfA}$
  factors as
  \begin{equation*}
    M' \xra{\psi'} M \xra{\psi} M''
  \end{equation*}
  where $M$ is in $\capUU$, then
  $\fuTdot(\grf) = \fuTdot(\psi\psi') = \fuTdot(\psi)\fuTdot(\psi')$
  factors through the contractible complex $\fuTdot(M)$, so one has
  $\fuTdot(\grf) = [0][0] = [0]$.
\end{prf*}

\begin{rmk}
  \label{rmk:lift}
  Let $\sfU$ be an additively closed subcategory of $\sfA$.

  Let $M$ be a right $\sfU$-Gorenstein object in $\sfA$ and $T$ a
  right $\sfU$-totally acyclic complex with $\Cy[0]{T}\is M$. It
  follows from \prpref{lift} that $T$ and $\fuTdot(M)$ are isomorphic
  in $\catKRtac{\capUU}$.

  Let $\mapdef{\grf}{M}{M'}$ be a morphism of right $\sfU$-Gorenstein
  objects in $\sfA$. For every morphism
  $\mapdef{\phi}{\fuTdot(M)}{\fuTdot(M')}$ that lifts $\grf$,
  \prpref{lift} yields $[\phi] = \fuTdot(\grf)$.
\end{rmk}

Let $\sfU$ be an additively closed subcategory of $\sfA$.  Recall from
\thmref{Frobenius-R} that $\catGR{\sfA}$ is a Frobenius category with
$\capUU$ the subcategory of projective-injective objects. Denote by
$\stabcatGR{\sfA}$ the associated stable category. It is a
triangulated category, see for example Krause~\cite[7.4]{HKr07a}, and
it is immediate from \prpref{Tdot} that $\fuTdot$ induces a
triangulated functor
$\mapdef[\lra]{\fuT}{\stabcatGR{\sfA}}{\catKRtac{\capUU}}$.

\begin{thm}
  \label{thm:T}
  Let $\sfU$ be an additively closed subcategory of $\sfA$.  There is
  a biadjoint triangulated equivalence
  \begin{equation*}
    \xymatrix@C=5.1pc{
      \stabcatGR{\sfA} \ar@<3.5pt>[r]^-{\fuT} &
      \catKRtac{\capUU} \ar@<3.5pt>[l]^-{\mathrm{Z}_0}\:.
    }
  \end{equation*}
\end{thm}

\begin{prf*}
  Set $\sfW = \capUU$ and notice that $\sfW$ is additively closed. The
  functors $\fuT$ and $\mathrm{Z}_0$ are triangulated. We prove that
  $(\fuT,\mathrm{Z}_0)$ is an adjoint pair; a parallel argument shows
  that $(\mathrm{Z}_0,\fuT)$ is an adjoint pair. Let $M$ be a right
  $\sfU$-Gorenstein object and $T$ be a right $\sfU$-totally acyclic
  complex. The assignment $[\phi] \longmapsto [\Cy[0]{\phi}]$ defines
  a map
  \begin{equation*}
    \dmapdef{\Phi^{M,T}}{\cathom{\catKRtac{\sfW}}{\fuT(M)}{T}}
    {\cathom{\stabcatGR{\sfA}}{M}{\Cy[0]{T}}}\:.
  \end{equation*}
  By \lemref{lift0} there is a morphism of $\sfA$-complexes
  $\mapdef{\gre^T}{\fuT(\Cy[0]{T})}{T}$.  The assignment
  $[\grf] \longmapsto [\gre^T]\fuT(\grf)$ defines a map $\Psi^{M,T}$
  in the opposite direction.

  Let $[\phi] \in \cathom{\catKRtac{\sfW}}{\fuT(M)}{T}$ be given. Let
  $\mapdef{\phi_{M,T}}{\fuT(M)}{\fuT(\Cy[0]{T})}$ be a representative
  of the homotopy class $\fuT(\Cy[0]{\phi})$, cf.~\rmkref{lift}. The
  composite $\gre^T\phi_{M,T}$ agrees with $\phi$ on
  $M = \Cy[0]{\fuT(M)}$, so \lemref{lift} yields
  \begin{equation*}
    [\phi] = [\gre^T\phi_{M,T}] = [\gre^T]\fuT(\Cy[0]{\phi})
    = \Psi^{M,T}\Phi^{M,T}([\phi])\:.
  \end{equation*}
  Now let $[\grf] \in \cathom{\stabcatGR{\sfA}}{M}{\Cy[0]{T}}$ be
  given. Let $\mapdef{\grf_{M,T}}{\fuT(M)}{\fuT(\Cy[0]{T})}$ be a lift
  of $\grf$; that is, a representative of the homotopy class
  $\fuT(\grf)$. One now has
  \begin{align*}
    \Phi^{M,T}\Psi^{M,T}([\grf]) &= \Phi^{M,T}([\gre^T]\fuT(\grf)) \\
                                 & = \Phi^{M,T}([\gre^T\grf_{M,T}]) \\
                                 & = [\Cy[0]{\gre^T\grf_{M,T}}] \\
                                 & = [\Cy[0]{\gre^T}\Cy[0]{\grf_{M,T}}] \\
                                 & = [1^{\Cy[0]{T}}\grf] \\
                                 & = [\grf]\:.
  \end{align*}
  Thus $\Phi^{M,T}$ is an isomorphism.

  The unit of the adjunction is the identity as one has
  $\Cy[0]{\fuT(-)} = 1^{\stabcatGR{\sfA}}$, and it is straightforward
  to check that $\gre^T$ defined above determines the counit
  $\mapdef{\gre}{\fuT(\Cy[0]{-})}{1^{\catKRtac{\sfW}}}$.  To show that
  $\gre$ is an isomorphism, let $T\in\catKRtac{\sfW}$ be given and
  consider a lift of the identity
  $\Cy[0]{T} \to \Cy[0]{\fuT(\Cy[0]{T})}$ to a morphism
  $\mapdef{\gri^T}{T}{\fuT(\Cy[0]{T})}$; see \lemref{lift0}. The
  composite $\gre^T\gri^T$ agrees with $1^T$ on $\Cy[0]{T}$, so
  $\gre^T\gri^T$ is a homotopy equivalence by
  \lemref{lift}. Similarly, $\gri^T\gre^T$ is a homotopy
  equivalence. It follows that $\gre^T$ is a homotopy equivalence,
  i.e.\ $[\gre^T]$ is an isomorphism in $\catKRtac{\sfW}$.
\end{prf*}

\begin{cor}
  \label{cor:s-o}
  Let $\sfUV$ be a cotorsion pair in $\sfA$.  There is a biadjoint
  triangulated equivalence
  \begin{equation*}
    \xymatrix@C=5.1pc{
      \stabcatG[\capUV]{\sfA} \ar@<3.5pt>[r]^-{\fuT} &
      \catKtac{\capUV} \ar@<3.5pt>[l]^-{\mathrm{Z}_0}\:.
    }
  \end{equation*}
\end{cor}

\begin{prf*}
  This is \thmref{T} applied to the self-orthogonal additively closed
  subcategory $\capUV$ and written in the notation from
  \dfnref[Definitions~]{Gor} and \dfnref[]{Ktac}.
\end{prf*}

\begin{exa}
  \label{exa:std}
  Applied to the cotorsion pair $(\sfA, \catInj{\sfA})$, \corref{s-o}
  recovers the well-known equivalence of the stable category of
  Gorenstein injective objects and the homotopy category of totally
  acyclic complexes of injective objects; see
  \prpcite[7.2]{HKr05}. Applied to the cotorsion pair
  $(\catPrj{\sfA},\sfA)$, the corollary yields the corresponding
  equivalence
  $\stabcatG[\mathsf{Prj}]{\sfA} \qis \catKtac{\catPrj{\sfA}}$.
\end{exa}

\begin{rmk}
  Let $\sfU$ and $\sfV$ be additively closed subcategories of
  $\sfA$. In \lemref[]{lift0}--\thmref[]{T} we have focused on right
  $\sfU$-totally acyclic complexes and right $\sfU$-Gorenstein
  objects. There are, of course, parallel results about left
  $\sfV$-totally acyclic complexes and left $\sfV$-Gorenstein
  objects. In particular, there is a biadjoint triangulated
  equivalence
  \begin{equation*}
    \xymatrix@C=5.1pc{
      \stabcatGL{\sfA} \ar@<3.5pt>[r]^-{\operatorname{T_L}} &
      \catKLtac{\capVV} \ar@<3.5pt>[l]^-{\mathrm{Z}_0}\:.
    }
  \end{equation*}
  Notice that applied to a cotorsion pair $\sfUV$ this also yields
  \corref{s-o}.
\end{rmk}

%%%%%%%%%%%%%%%%%%%%%%%%%%%%%%%%%%%%%%%%%%%%%%%%%%%%%%%%%%%%%%%%%%%%%%
\section{Gorenstein flat-cotorsion modules}

\noindent
In this section and the next, $A$ is an associative ring. We adopt the
convention that an $A$-module is a left $A$-module; right $A$-modules
are considered to be modules over the opposite ring $\Aop$. The
category of $A$-modules is denoted $\catMod{A}$.

Given a cotorsion pair $\sfUV$ in $\catMod{A}$ the natural categories
of Gorenstein objects to consider are $\catGR{A}$, $\catGL{A}$, and
$\catG[(\capUV)]{A}$; see \rmkref{UVW}. For each of the absolute
cotorsion pairs $(\catPrj{A},\catMod{A})$ and
$(\catMod{A},\catInj{A})$, two of these categories of Gorenstein
objects coincide and the third is trivial. We start this section by
recording the non-trivial fact that the cotorsion pair
$(\catFlat{A},\catCot{A})$ exhibits the same behavior. For brevity we
denote the self-orthogonal subcategory $\catFlat{A}\cap\catCot{A}$ of
flat-cotorsion modules by $\catFlatCot{A}$.

Bazzoni, Cort\'es Izurdiaga, and Estrada \thmcite[4.1]{BCE} prove:

\begin{fct}
  \label{cot-cyc}
  An acyclic complex of cotorsion $A$-modules has cotorsion cycle
  modules.
\end{fct}

\begin{prp}
  \label{prp:flatcot}
  A $\catFlatCot{A}$-totally acyclic complex is right
  $\catFlat{A}$-totally acyclic, and a left $\catCot{A}$-totally acyclic
  complex is contractible. In particular, one has
  \begin{equation*}
    \catGR[\sfFlat]{A} = \catG[\sfFlatCot]{A} \qqand 
    \catGL[\sfCot]{A} = \catFlatCot{A}\;.
  \end{equation*}
\end{prp}

\begin{prf*}
  In a $\catFlatCot{A}$-totally acyclic complex, the cycle modules are
  cotorsion by \ref{cot-cyc}, whence the complex is right
  $\catFlat{A}$-totally acyclic by \dfnref{tac}.  By \rmkref{UVW}
  every right $\catFlat{A}$-totally acyclic complex is
  $\catFlatCot{A}$-totally acyclic, so the first equality of
   categories follows from \dfnref{Gor}. In a left
  $\catCot{A}$-totally acyclic complex, the cycle modules are
  flat-cotorsion, again by \dfnref{tac} and \ref{cot-cyc}, so such a
  complex is contractible, and the second equality follows.
\end{prf*}

We introduce a less symbol-heavy terminology.

\begin{dfn}
  \label{dfn:FC-tac}
  A $\catFlatCot{A}$-totally acyclic complex is called a \emph{totally
    acyclic complex of flat-cotorsion modules.}  A cycle module in
  such a complex, that is, a $\catFlatCot{A}$-Gorenstein module, is
  called a \emph{Gorenstein flat-cotorsion module.}
\end{dfn}

Recall that a complex $T$ of flat $A$-modules is called \emph{\Ftac}
if it is acyclic and the complex $\tp{I}{T}$ is acyclic for every
injective $\Aop$-module $I$.

\begin{thm}
  \label{thm:2Ftac}
  Let $A$ be right coherent. For an $A$-complex $T$ the following
  conditions are equivalent
  \begin{eqc}
  \item $T$ is a totally acyclic complex of flat-cotorsion modules.
  \item $T$ is a complex of flat-cotorsion modules and \Ftac.
  \item $T$ is right $\catFlat{A}$-totally acyclic.
  \end{eqc}
%  In particular, one has $\catGR[\sfFlat]{A} = \catG[\sfFlatCot]{A}$.
\end{thm}

\begin{prf*}
  Per \rmkref{UVW} condition \eqclbl{iii} implies \eqclbl{i}.

  \proofofimp{i}{ii} If $T$ is a totally acyclic complex of
  flat-cotorsion modules, then by \prpref{tac} it is an acyclic
  complex of flat-cotorsion modules. For every injective $\Aop$-module
  $I$ the $A$-module $\Hom[\ZZ]{I}{\QQ/\ZZ}$ is flat-cotorsion, as $A$
  is right coherent. Now it follows by the isomorphism
  \begin{equation*}
    \tag{$\ast$}
    \Hom[A]{T}{\Hom[\ZZ]{I}{\QQ/\ZZ}} \is \Hom[\ZZ]{\tp[A]{I}{T}}{\QQ/\ZZ}
  \end{equation*}
  and faithful injectivity of $\QQ/\ZZ$ that $\tp{I}{T}$ is acyclic.

  \proofofimp{ii}{iii} If $T$ is a complex of flat-cotorsion
  $A$-modules and \Ftac, then $T$ satisfies conditions $(\scR.1)$ and
  $(\scR.2)$ in \dfnref{tac}. By \ref{cot-cyc} the cycles modules of
  $T$ are cotorsion, so $T$ also satisfies condition
  $(\scR.3)$. Further, as $A$ is right coherent, every flat-cotorsion
  $A$-module is a direct summand of a module of the form
  $\Hom[\ZZ]{I}{\QQ/\ZZ}$, for some injective $\Aop$-module $I$; see
  e.g.\ Xu \lemcite[3.2.3]{xu}. Now it follows from the isomorphism
  $(\ast)$ that $\Hom{T}{W}$ is acyclic for every
  $W\in \catFlatCot{A}$. That is, $T$ also satisfies condition
  \dfnref[]{tac}(\scR.4).
%
%  The last assertion now follows from \dfnref{Gor} and
%  \rmkref{R=LGor}.
\end{prf*}

Recall that an $A$-module $M$ is called \emph{Gorenstein flat} if
there exists an \Ftac\ complex $F$ of flat $A$-modules with
$\Cy[0]{F} = M$.  The full subcategory of $\catMod{A}$ whose objects
are the Gorenstein flat modules is denoted $\catGFlat{A}$.

Gillespie \corcite[3.4]{JGl17} proved that the category
$\catCot{A} \cap \catGFlat{A}$ is Frobenius if $A$ is right
coherent. That it remains true without the coherence assumption is an
immediate consequence of \corcite[3.12]{JSrJSt} discussed
\emph{ibid.}; for convenience we include the statement as part of the
next result.

\begin{thm}
  \label{thm:perfect}
  The category $\catCot{A} \cap \catGFlat{A}$ is Frobenius and
  $\catFlatCot{A}$ is the subcategory of projective-injective
  objects. Moreover, the following conditions are equivalent.
  \begin{eqc}
  \item $A$ is left perfect.
  \item The category $\catGFlat{A}$ is Frobenius.
  \item One has $\catGFlat{A} = \catCot{A} \cap \catGFlat{A}$.
  \end{eqc}
  Furthermore, if $A$ is right coherent then these conditions are
  equivalent to
  \begin{eqc}\setcounter{eqc}{3}
  \item An $A$-module is Gorenstein flat if and only if it is
    Gorenstein projective.
  \end{eqc}
\end{thm}

\begin{prf*}
  By \corcite[3.12]{JSrJSt} the category $\catGFlat{A}$ is closed
  under extensions, and $\catGFlat{A} \cap \catGFlat{A}^\perp$ is the
  subcategory $\catFlatCot{A}$ of flat-cotorsion modules. It follows
  that $\catCot{A} \cap \catGFlat{A}$ is closed under extensions, and
  that modules in $\catFlatCot{A}$ are both projective and injective
  in $\catCot{A} \cap \catGFlat{A}$. Let $P$ be a projective object in
  $\catCot{A} \cap \catGFlat{A}$; it fits in an exact sequence
  \begin{equation*}
    \tag{$\ast$}
    0 \lra C \lra F \lra P \lra 0
  \end{equation*}
  where $F$ is flat and $C$ is cotorsion; see Bican, El Bashir, and
  Enochs \cite{BEE-01}. As $P$ is cotorsion it follows that $F$ is
  flat-cotorsion. By \corcite[3.12]{JSrJSt} the category
  $\catGFlat{A}$ is resolving, so $C$ is Gorenstein flat. Thus,
  $(\ast)$ is an exact sequence in $\catCot{A} \cap \catGFlat{A}$,
  whence it splits by the assumption on $P$. In particular, $P$ is
  flat-cotorsion. Now let $I$ be an injective object in
  $\catCot{A} \cap \catGFlat{A}$. It fits by \corcite[3.12]{JSrJSt} in
  an exact sequence
  \begin{equation*}
    \tag{$\dagger$}
    0 \lra I \lra F \lra G \lra 0
  \end{equation*}
  where $F$ belongs to $\catGFlat{A}^\perp$ and $G$ is Gorenstein
  flat. It follows that $F$ is Gorenstein flat and hence
  flat-cotorsion, still by \corcite[3.12]{JSrJSt}. Finally, $G$ is
  cotorsion as both $I$ and $F$ are cotorsion. Thus, $(\dagger)$ is an
  exact sequence in $\catCot{A} \cap \catGFlat{A}$, whence it splits
  by the assumption on $I$. In particular, $I$ is flat-cotorsion.

  \proofofimp{i}{iii} Assuming that $A$ is left perfect, every flat
  $A$-module module is projective, whence every $A$-module is
  cotorsion.

  \proofofimp{iii}{ii} Evident as $\catCot{A} \cap \catGFlat{A}$ is
  Frobenius as shown above.

  \proofofimp{ii}{i} Assume that $\catGFlat{A}$ is Frobenius and
  denote by $\sfW$ its subcategory of projective-injective objects. To
  prove that $A$ is left perfect it suffices by a result of Guil
  Asensio and Herzog \corcite[20]{PGAIHr05} to show that the free
  module $A^{(\NN)}$ is cotorsion.  As $A^{(\NN)}$ is flat, in
  particular Gorenstein flat, and as $\catGFlat{A}$ by assumption has
  enough projectives, there is an exact sequence
  $0 \to K \to W \to A^{(\NN)} \to 0$ with $W\in \sfW$. The sequence
  splits because $A^{(\NN)}$ is projective, so it suffices to show
  that modules in $\sfW$ are cotorsion. Fix $W\in\sfW$, let $F$ be a
  flat $A$-module, and consider an extension
  \begin{equation*}
    \tag{$\ddagger$}
    0 \lra W \lra E \lra F \lra 0\:.
  \end{equation*}
  As $\catGFlat{A}$ by \corcite[3.12]{JSrJSt} is closed under
  extensions, the module $E$ is Gorenstein flat. As $W$ is injective
  in $\catGFlat{A}$ it follows that the sequence $(\ddagger)$ splits,
  i.e.\ one has $\Ext{1}{F}{W}=0$. That it, $W$ is cotorsion.

  \proofofimp{iv}{ii} By \thmref{Frobenius-R} the category of
  Gorenstein projective $A$-modules is Frobenius.

  \proofofimp{i}{iv} If $A$ is perfect and right coherent, then it
  follows from \thmref{2Ftac} that an $A$-module is Gorenstein flat if
  and only if it is Gorenstein projective.
\end{prf*}

By \thmref{perfect} the category $\catGFlat{A}$ is only Frobenius when
every $A$-module is cotorsion, and the take-away is that the
appropriate Frobenius category to focus on is
$\catCot{A} \cap \catGFlat{A}$.  If $A$ is right coherent ring, then
this category contains $\catG[\sfFlatCot]{A}$, by \thmref{2Ftac} and
\pgref{cot-cyc}, and one goal of the next section is to prove the
reverse inclusion; that is \thmref{GG}.

%%%%%%%%%%%%%%%%%%%%%%%%%%%%%%%%%%%%%%%%%%%%%%%%%%%%%%%%%%%%%%%%%%%%%%
\section{The stable category of Gorenstein flat-cotorsion modules}

\noindent
Recall that an $A$-complex $P$ is called \emph{pure-acyclic} if the
complex $\tp[A]{N}{P}$ is acyclic for every $\Aop$-module $N$. In
particular, an acyclic complex $P$ of flat $A$-modules is pure-acyclic
if and only if all cycle modules $\Cy[i]{P}$ are flat.

\begin{fct}
  \label{cot-approx}
  Let $M$ be an $A$-complex. It follows\footnote{Although
    \cite[cor. 4.10]{JGl04} is stated for commutative rings, it is
    standard that the result remains valid
    without this assumption; see for example the discussion before
    \thmcite[4.2]{SEsJGl19}.} from Gillespie~\corcite[4.10]{JGl04} that
  there exists an exact sequence of $A$-complexes
  \begin{equation*}
    0 \lra M \lra C \lra P \lra 0
  \end{equation*}
  where $C$ is a complex of cotorsion modules and $P$ is a
  pure-acyclic complex of flat modules.
\end{fct}

The first theorem of this section shows that if $A$ is right coherent,
then the cotorsion modules in $\catGFlat{A}$ are precisely the
non-trivial Gorenstein modules associated to the cotorsion pair
$(\catFlat{A},\catCot{A})$; namely the Gorenstein flat-cotorsion
modules or, equivalently, the right $\catFlat{A}$-Gorenstein modules.

\begin{thm}
  \label{thm:GG}
  Let $A$ be right coherent. There are equalities
  \begin{equation*}
    \catCot{A} \cap \catGFlat{A}
    =\catG[\sfFlatCot]{A}
    =\catGR[\sfFlat]{A}\:.
  \end{equation*}
\end{thm}

\begin{prf*}
  The second equality is by \prpref{flatcot}, and the containment
  \begin{equation*}
    \catCot{A} \cap \catGFlat{A} \supseteq \catG[\sfFlatCot]{A}
  \end{equation*}
  follows from \pgref{cot-cyc} and \thmref{2Ftac}. It remains to show
  the reverse containment.

  Let $M$ be a Gorenstein flat $A$-module that is also cotorsion. By
  definition, there is an \Ftac\ complex $F$ of flat $A$-modules with
  $\Cy[0]{F} = M$. Further, \pgref{cot-approx} yields an exact
  sequence of $A$-complexes
  \begin{equation*}
    \tag{1}
    0 \lra F \xra{\gri'} T \xra{\pi'} P \lra 0
  \end{equation*}
  where $T$ is a complex of cotorsion modules and $P$ is a
  pure-acyclic complex of flat modules. It follows that $T$ is a
  complex of flat modules; moreover, since $P$ is trivially \Ftac, so
  is $T$. As $A$ is right coherent, it now follows from \thmref{2Ftac}
  that $T$ is a totally acyclic complex of flat-cotorsion modules.

  The functor $\Cy{-}$ is left exact, and since $F$ is acyclic a
  standard application of the Snake Lemma yields the exact sequence
  \begin{equation*}
    \tag{2}
    0 \lra M \xra{\gri} \Cy[0]{T} \xra{\pi} \Cy[0]{P} \lra 0
  \end{equation*}
  where $\gri$ and $\pi$ are the restrictions of the morphisms from
  $(1)$. As $M$ is cotorsion and $\Cy[0]{P}$ is flat, $(2)$
  splits. Set $Z=\Cy[0]{P}$ and denote by $\varrho$ the section with
  $\pi\varrho = 1^Z$.  By \pgref{cot-cyc} the module $\Cy[0]{T}$ is
  cotorsion, so it follows that $Z$ is a flat-cotorsion module. Now,
  as $\Cy[-1]{P}$ is flat, the exact sequence
  \begin{equation*}
    0 \lra Z \xra{\gre^P_0} P_0 \lra \Cy[-1]{P} \lra 0
  \end{equation*}
  splits; denote by $\grs$ the section with $\grs\gre_0^P = 1^Z$. By
  commutativity of the diagram
  \begin{equation*}
    \xymatrix{
      \Cy[0]{T} \ar[r]^-{\pi} \ar[d]^{\gre_0^T}
      & Z \ar[d]^{\gre^P_0}\\
      T_0 \ar[r]^{\pi'_0}  & P_0
    }
  \end{equation*}
  one has $\grs\pi'_0\gre_0^T\varrho = \grs\gre_0^P\pi\varrho =
  1^Z$. It follows that $\mapdef{\grs\pi'_0}{T_0}{Z}$ is a split
  surjection with section $\gre_0^T\varrho$.

  As $Z$ is flat and $\Cy[1]{T}$ is cotorsion, there is an exact
  sequence
  \begin{equation*}
    0 \lra \Hom{Z}{\Cy[1]{T}} \lra \Hom{Z}{T_1}
    \lra \Hom{Z}{\Cy[0]{T}} \lra 0\:.
  \end{equation*}
  It follows that there is a homomorphism $\mapdef{\grz}{Z}{T_1}$ with
  $\partial_1^T\grz = \varrho$ and, therefore,
  $\partial_1^T\grz = \gre_0^T\varrho$ as homomorphisms from $Z$ to
  $T_0$. As $\partial_0^T\gre_0^T\varrho=0$ trivially holds, the
  homomorphisms $\grz$ and $\gre_0^T\varrho$ yield a morphism of
  complexes:
  \begin{equation*}
    \xymatrix{
      D \ar[d]^-\rho & \hspace{-2pc}= \ \cdots \ \ar[r] & 0 \ar[r] \ar[d]
      & Z \ar[r]^-{=} \ar[d]^-{\grz}
      & Z \ar[r] \ar[d]^-{\gre_0^T\varrho} & 0\ar[r] \ar[d] & \ \cdots \\
      T & \hspace{-2pc}= \ \cdots \ar[r] & T_2 \ar[r]^{\partial_2^T} & T_1
      \ar[r]^{\partial_1^T} & T_0 \ar[r]^{\partial_0^T}
      & T_{-1}\ar[r] & \cdots }
  \end{equation*}
  This is evidently a split embedding in the category of complexes
  whose section given by the homomorphisms
  $\mapdef{\grs\pi'_0\partial_1^T}{T_1}{Z}$ and
  $\mapdef{\grs\pi'_0}{T_0}{Z}$. The restriction of the split exact
  sequence of complexes
  \begin{equation*}
    \tag{3}
    0 \lra D \lra T \lra T' \lra 0
  \end{equation*}
  to cycles is isomorphic to the split exact sequence
  $0 \lra Z \xra{\varrho} \Cy[0]{T} \lra M \lra 0$, see $(2)$, so it
  follows that the complex $T'$ has $\Cy[0]{T'} \is M$.

  In $(3)$ both $D$ and $T$ are complexes of flat-cotorsion modules
  and \Ftac, so also $T'$ is a complex of flat-cotorsion modules and
  \Ftac. Now it follows from \thmref{2Ftac} that $T'$ is a totally
  acyclic complex of flat-cotorsion modules, whence the module
  $M \is \Cy[0]{T'}$ is Gorenstein flat-cotorsion.
\end{prf*}

\begin{cor}
  \label{cor:eq0}
  Let $A$ be right coherent. There is a triangulated equivalence
  \begin{equation*}
    \stabcatG[\sfFlatCot]{A} \qis \catKFtac{\catFlatCot{A}}\:.
  \end{equation*}
\end{cor}

\begin{prf*}
  Immediate from \thmref[Theorems~]{T}, \thmref[]{2Ftac}, and
  \thmref[]{GG}; see also the diagram in \pgref{summary}.
\end{prf*}

\begin{cor}
  Let $A$ be right coherent. The category $\catG[\sfFlatCot]{A}$ is
  closed under direct summands.
\end{cor}

\begin{prf*}
  Immediate from the theorem as both $\catCot{A}$ and $\catGFlat{A}$
  are closed under direct summands; for the latter
  see~\corcite[3.12]{JSrJSt}.
\end{prf*}

Gorenstein flat $A$-modules are, within the framework of Sections
1--2, not born out of a cotorsion pair, not even out of a
self-orthogonal subcategory of $\catMod{A}$. However, they form the
left half of a cotorsion pair, and also out of that pair comes the
Gorenstein flat-cotorsion modules.

\begin{rmk}
  Let $A$ be right coherent. Enochs, Jenda, and Lopez-Ramos
  \thmcite[2.11]{EJL-04} show that $\catGFlat{A}$ is the left half of
  a cotorsion pair, and Gillespie \prpcite[3.2]{JGl17} shows that
  $\catGFlat{A} \cap \catGFlat{A}^\perp$ is
  $\catFlatCot{A}$\footnote{\,\v{S}aroch and
    {\v{S}}tov{\'{\i}}{\v{c}}ek~\corcite[3.12]{JSrJSt} show that all
    of this is true without assumptions on $A$, and we used that
    crucially in the proof of \thmref{perfect}. The results from
    \cite{EJL-04} and \cite{JGl17} suffice to prove \thmref[]{perfect}
    for a right coherent ring.}

  A right $\catGFlat{A}$-totally acyclic complex as well as a left
  $\catGFlat{A}^\perp$-totally acyclic complex is by \rmkref{UVW} and
  \dfnref{FC-tac} a totally acyclic complex of flat-cotorsion
  modules. For a right $\catGFlat{A}$-totally acyclic complex $T$, it
  follows from \dfnref{tac} that $\Hom{G}{T}$ is acyclic for every
  Gorenstein flat $A$-module $G$, in particular for every Gorenstein
  flat-cotorsion module. That is, such a complex is contractible and,
  therefore, a right $\catGFlat{A}$-Gorenstein module is
  flat-cotorsion. On the other hand, the cycles in a left
  $\catGFlat{A}^\perp$-totally acyclic complex are by \dfnref{tac}
  Gorenstein flat and by \pgref{cot-cyc} cotorsion, so a left
  $\catGFlat{A}^\perp$-Gorenstein module is by \thmref{GG} Gorenstein
  flat-cotorsion.
\end{rmk}

Let $\catKpac{\catFlat{A}}$ denote the full subcategory of
$\catK{\catFlat{A}}$ whose objects are pure-acyclic; notice that it is
contained in $\catKFtac{\catFlat{A}}$.  Via \pgref{cot-cyc} and the
dual of \pgref{cot-approx} one could obtain the next theorem as a
consequence of a standard result \prpcite[10.2.7]{cas}; we opt for a
direct argument.

\begin{thm}
  \label{thm:I}
  The composite
  \begin{equation*}
    \fuI\colon
    \catKFtac{\catFlatCot{A}} \lra
    \catKFtac{\catFlat{A}} \lra
    \frac{\catKFtac{\catFlat{A}}}{\catKpac{\catFlat{A}}}
  \end{equation*}
  of canonical functors is a triangulated equivalence of categories.
\end{thm}

\begin{prf*}
  Let $\fuI$ be the composite of the inclusion followed by Verdier
  localization; notice that $\fuI$ is the identity on objects. We
  argue that the functor $\fuI$ is essentially surjective, full, and
  faithful.

  Let $F$ be an \Ftac\ complex of flat modules. By \ref{cot-approx}
  there is an exact sequence
  \begin{equation*}
    \tag{$\ast$}
    0 \lra F \lra C^F \lra P^F \lra 0
  \end{equation*}
  where $C^F$ is a complex of cotorsion modules and $P^F$ is in
  $\catKpac{\catFlat{A}}$.  As $F$ and $P^F$ are \Ftac\ complexes of
  flat $A$-modules so is $C^F$; that is, $C^F$ belongs to
  $\catKFtac{\catFlatCot{A}}$. It follows from ($\ast$) that $F$ and
  $C^F$ are isomorphic in the Verdier quotient
  $\frac{\catKFtac{\catFlat{A}}}{\catKpac{\catFlat{A}}}$. Thus $\fuI$
  is essentially surjective.

  Let $F$ and $F'$ be \Ftac\ complexes of flat-cotorsion modules. A
  morphism $F \to F'$ in
  $\frac{\catKFtac{\catFlat{A}}}{\catKpac{\catFlat{A}}}$ is a diagram
  in $\catKFtac{\catFlat{A}}$
  \begin{equation*}
    \tag{$\ast$}
    F  \xrightarrow{\ [\gra] \ } X \xleftarrow[\qis]{\ [\grf] \ } F'
  \end{equation*}
  such that the complex $\Cone{\grf}$ belongs to
  $\catKpac{\catFlat{A}}$. Let $\gri$ be the embedding $X \to C^X$
  from \ref{cot-approx}. It is elementary to verify that the composite
  $\mapdef{\gri\grf}{F'}{C^X}$ has a pure-acyclic mapping cone; see
  \lemcite[2.7]{LWCPTh}. Since $F'$ and $C^X$ are complexes of
  flat-cotorsion modules, so is $\Cone{\gri\grf}$. It now follows by
  way of \ref{cot-cyc} that $\Cone{\gri\grf}$ is contractible; that
  is, $\gri\grf$ is a homotopy equivalence. Thus $[\gri\grf]$ has an
  inverse in $\catKFtac{\catFlat{A}}$, i.e.\
  $[\gri\grf]^{-1} = [\psi]$ for some morphism
  $\mapdef{\psi}{C^X}{F'}$. The commutative diagram
  \begin{equation*}
    \begin{gathered}
      \xymatrix@=2.8pc{
        {} & X \ar[d]^-{[\gri]}_-\qis  & {} \\
        F \ar[r]^-{[\gri\gra]} \ar[ur]^-{[\gra]}
        \ar[dr]_-{[\psi\gri\gra]}^-\qis & C^X & F'
        \ar[l]_-{[\gri\grf]}^-\qis \ar[ul]_-{[\grf]}^-\qis
        \ar[dl]^-{[1^{F'}]}_-\qis \\
        {} & F' \ar[u]_-{[\gri\grf]}^-\qis & {} }
    \end{gathered}
  \end{equation*}
  now shows that the morphism $(\ast)$ is equivalent to
  $F \xra{[\psi\gri\gra]} F' \xla{[1^{F'}]} F'$, which is
  $\fuI(\psi\gri\gra)$. This shows that $\fuI$ is full.

  Finally, let $\mapdef{\gra}{F}{F'}$ be a morphism of \Ftac\
  complexes of flat-cotorsion modules, and assume that $\fuI([\gra])$
  is zero. It follows that there is a commutative diagram in
  $\catKFtac{\catFlat{A}}$,
  \begin{equation*}
    \begin{gathered}
      \xymatrix@=2.8pc{
        {} & F' \ar[d]^-{[\grf]}_-\qis  & {} \\
        F \ar[r]^-{[\grf\gra]} \ar[ur]^-{[\gra]} \ar[dr]_-{[0]}^-\qis
        & X & F' \ar[l]_-{[\grf]}^-\qis \ar[ul]_-{[1^{F'}]}^-\qis
        \ar[dl]^-{[1^{F'}]}_-\qis \\
        {} & F' \ar[u]_-{[\grf]}^-\qis & {} }
    \end{gathered}
  \end{equation*}
  where the mapping cone of $\grf$ is in $\catKpac{\catFlat{A}}$. The
  diagram yields $[\grf\gra]=[0]$ and, therefore,
  $[\gri\grf][\gra] = [\gri\grf\gra] = [0]$ where $\gri$ is the
  embedding $X\to C^X$ from \pgref{cot-approx}. As above, $[\gri\grf]$
  is invertible in $\catKFtac{\catFlat{A}}$, so one has $[\gra]=[0]$
  in $\catKFtac{\catFlat{A}}$. That is, $\gra$ is null-homotopic, and
  hence $[\gra]=0$ in $\catKFtac{\catFlatCot{A}}$.
\end{prf*}

\begin{bfhpg}[Summary]
  \label{summary}
  Let $A$ be right coherent. By \thmref[Theorems~]{T} and \thmref[]{I}
  there are triangulated equivalences
  \begin{equation*}
    \xymatrix@C=2.5pc@R=1.5pc{
      & \catKFtac{\catFlatCot{A}}\ar@{=}[d] \ar[r]^-\fuI_-\qis
      & \smash{\displaystyle\frac{\catKFtac{\catFlat{A}}}{\catKpac{\catFlat{A}}}}\\
      \stabcatGR[\sfFlat]{A} \ar[r]^-{\fuT}_-\qis
      & \catKRtac[\sfFlat]{\catFlatCot{A}} \\
      \stabcatG[\sfFlatCot]{A}  \ar@{=}[u]
      & \catKtac{\catFlatCot{A}} \ar@{=}[u]
    }
  \end{equation*}
  where the equalities come from \prpref{flatcot} and
  \thmref[Theorems~]{2Ftac} and \thmref[]{GG}.
\end{bfhpg}

\begin{cor}
  \label{cor:eq}
  Let $A$ be right coherent. There is a triangulated equivalence
  \begin{equation*}
    \stabcatG[\sfFlatCot]{A} \qis
    \frac{\catKFtac{\catFlat{A}}}{\catKpac{\catFlat{A}}}\:.
  \end{equation*}
\end{cor}

\begin{prf*}
  See the diagram in \pgref{summary}.
\end{prf*}

In the special case where $A$ is commutative noetherian of finite
Krull dimension, the next result is immediate from
\lemcite[4.22]{DMfSSl11} and \corref{eq}.

\begin{cor}
  \label{cor:GPGF}
  Let $A$ be right coherent ring such that all flat $A$-modules have
  finite projective dimension. There is a triangulated equivalence of
  categories
  \begin{equation*}
    \stabcatG[\mathsf{Prj}]{A} \qis \stabcatG[\sfFlatCot]{A}\:.
  \end{equation*}
\end{cor}

\begin{prf*}
  Under the assumptions on $A$, a complex of projective $A$-modules is
  totally acyclic if and only if it \Ftac; see \cite[claims 2.4 and
  2.5]{LWCKKt18}. By \thmcite[5.1]{SEsJGl19} there is now a
  triangulated equivalence of categories
  \begin{equation*}
    \catKtac{\catPrj{A}}\qis 
    \frac{\catKFtac{\catFlat{A}}}{\catKpac{\catFlat{A}}}\:.
  \end{equation*}
  Now apply the equivalences from \exaref{std} and \corref{eq}.
\end{prf*}

%%%%%%%%%%%%%%%%%%%%%%%%%%%%%%%%%%%%%%%%%%%%%%%%%%%%%%%%%%%%%%%%%%%%%%
\section*{Acknowledgments}

\noindent
We thank Petter Andreas Bergh, Tsutomu Nakamura, and Mark Walker for
conversations and comments on an early draft of this paper.

%%%%%%%%%%%%%%%%%%%%%%%%%%%%%%%%%%%%%%%%%%%%%%%%%%%%%%%%%%%%%%%%%%%%%%
% \bibliographystyle{amsplain-nodash} \bibliography{../+references}

\def\soft#1{\leavevmode\setbox0=\hbox{h}\dimen7=\ht0\advance \dimen7
  by-1ex\relax\if t#1\relax\rlap{\raise.6\dimen7
    \hbox{\kern.3ex\char'47}}#1\relax\else\if T#1\relax
  \rlap{\raise.5\dimen7\hbox{\kern1.3ex\char'47}}#1\relax \else\if
  d#1\relax\rlap{\raise.5\dimen7\hbox{\kern.9ex
      \char'47}}#1\relax\else\if D#1\relax\rlap{\raise.5\dimen7
    \hbox{\kern1.4ex\char'47}}#1\relax\else\if l#1\relax
  \rlap{\raise.5\dimen7\hbox{\kern.4ex\char'47}}#1\relax \else\if
  L#1\relax\rlap{\raise.5\dimen7\hbox{\kern.7ex
      \char'47}}#1\relax\else\message{accent \string\soft \space #1
    not defined!}#1\relax\fi\fi\fi\fi\fi\fi}
\providecommand{\MR}[1]{\mbox{\href{http://www.ams.org/mathscinet-getitem?mr=#1}{#1}}}
\renewcommand{\MR}[1]{\mbox{\href{http://www.ams.org/mathscinet-getitem?mr=#1}{#1}}}
\providecommand{\arxiv}[2][AC]{\mbox{\href{http://arxiv.org/abs/#2}{\sf
      arXiv:#2 [math.#1]}}} \def\cprime{$'$}
\providecommand{\bysame}{\leavevmode\hbox to3em{\hrulefill}\thinspace}
\providecommand{\MR}{\relax\ifhmode\unskip\space\fi MR }
\providecommand{\MRhref}[2]{%
  \href{http://www.ams.org/mathscinet-getitem?mr=#1}{#2} }
\providecommand{\href}[2]{#2}

\end{document}